\documentclass[10pt]{article}

\usepackage[superscript,sort]{cite}

\usepackage{ifthen,ifpdf}
\ifpdf
	\usepackage[pdftex]{hyperref}	 \hypersetup{colorlinks=true,linkcolor=black,citecolor=black,urlcolor=blue,backref=page,bookmarks=true,breaklinks=true,plainpages=false }
	\usepackage[pdftex]{graphicx}
	\usepackage[usenames,pdftex]{color}
\else
	\usepackage[colorlinks=true,breaklinks=true,plainpages=false]{hyperref}
	\usepackage{graphicx}
	\usepackage[usenames]{color}
\fi

\usepackage{amsthm,amscd,amsxtra,amsfonts,amsmath,amssymb,multirow}
\usepackage{wrapfig}
\usepackage[footnotesize]{caption}
\usepackage[tiny,compact]{titlesec}
\usepackage[textwidth=0.8in,textsize=footnotesize]{todonotes}
\usepackage{algorithm}
\usepackage[noend]{algpseudocode}

\usepackage{helvet}

\usepackage[letterpaper,top=0.5in,bottom=0.5in,left=0.5in,right=0.5in]{geometry}

\titlespacing*{\section}{0pt}{*0}{*0}
\titlespacing*{\subsection}{0pt}{*0}{*0}
\titlespacing*{\subsubsection}{0pt}{*0}{*0}
\titlespacing{\paragraph}{0pt}{*0}{*1}
\definecolor{MyPurple}{rgb}{1,0,1}

\newtheorem{remark}{Remark}

\makeatletter
\def\BState{\State\hskip-\ALG@thistlm}
\makeatother

\begin{document}

\pagenumbering{roman}

%

%

\clearpage \pagebreak \setcounter{page}{1}
\renewcommand{\thepage}{{\arabic{page}}}

\title{Parameter optimization in differential geometry based solvation models}

\author{
Bao Wang$^1$,    Nathan A. Baker$^2$,  and
G. W. Wei$^{1,3,4}$\footnote{Please address correspondence to Guowei Wei. E-mail:wei@math.msu.edu}
\\
$^1$Department of Mathematics, \\
Michigan State University, MI 48824, USA\\
$^2$Pacific Northwest National Laboratory, P.O. Box 999, MS K7-28, \\
Richland, WA 99352, USA \\
$^3$Department of Electrical and Computer Engineering,\\
Michigan State University, MI 48824, USA \\
$^4$Biochemistry  and Molecular Biology,\\
Michigan State University, MI 48824, USA
}

\date{\today}
\maketitle
\begin{abstract}

Differential geometry (DG)  based solvation models are a new class of variational implicit solvent approaches that are  able to avoid unphysical solvent-solute boundary definitions and associated geometric singularities, and dynamically couple polar and nonpolar interactions in a self-consistent framework. Our earlier study indicates that DG based nonpolar solvation model outperforms other methods in nonpolar solvation energy predictions. However, the DG based full solvation model has not shown its superiority in solvation analysis, due to its difficulty in parametrization, which must ensure the stability of the solution of strongly coupled nonlinear Laplace-Beltrami and Poisson-Boltzmann equations. In this work,  we introduce  new parameter learning algorithms based on perturbation and convex optimization theories to stabilize the numerical solution and thus achieve an optimal parametrization of the DG based solvation models.  An interesting feature of the present DG based solvation model is that it provides accurate solvation free energy predictions   for both polar and nonploar  molecules in a unified formulation. Extensive numerical experiment demonstrates that the present DG based solvation model delivers some of the most accurate predictions of the solvation free energies for a large number of molecules.

\end{abstract}

\newpage

\section{Introduction}

Biological processes, such as signaling, gene regulation,  transcription, translation, et cetera  govern the cell growth, cellular differentiation, fermentation, fertilization, germination, etc. in living organisms. Chemical processes, such as oxidation, reduction, hydrolysis, nitrification, polymerization, and so forth  underpin biological processes.  Physical processes, particularly solvation, are involved in all the aforementioned chemical and biological processes. Therefore, a prerequisite for the understanding of chemical and biological processes is to study the solvation process.  As a physical process, solvation does not involve the formation and/or breaking of any covalent bond, but is associated with solvent and solute electrostatic, dipolar, induced dipolar, and van der Waals interactions.

Experimentally, solvation can be analyzed by the measurement of solvation free energies. Theoretically, solavtion can be investigated by quantum mechanics, molecular mechanics, integral equation,  implicit solvent models,    and simple phenomenological modifications of Coulomb’s law. Among, the implicit solvent models are known to balance the computational complexity and the accuracy in the solvation free energy prediction, and thus, offer an efficient approach.

The general idea of  implicit solvent models is to treat the solvent as a dielectric continuum and describe the solute in atomistic detail \cite{Davis:1990a,Sharp:1990a,Honig:1995a,Roux:1999,Jinnouchi:2008}. The total solvation free energy is decomposed into nonpolar and polar parts. There is a wide variety of ways to carry out this decomposition. For example,  nonpolar energy contributions can be modeled in two stages:  the work of displacing solvent when adding a rigid solute to the solvent and the dispersive nonpolar interactions between the solute atoms and surrounding solvent. The polar part is due to the electrostatic interactions and can be approximated by generalized Born (GB) \cite{Dominy:1999,Bashford:2000,Tsui:2000,Onufriev:2002,Gallicchio:2002,Zhu:2005,Koehl:2006,Tjong:2007b,Mongan:2007,Chen:2008,Grant:2007}, polarizable continuum  (PC)  \cite{Tomasi:2005}
and Poisson-Boltzmann (PB) models \cite{Lamm:2003,Fogolari:2002,Sharp:1990a,Davis:1990a,Zhou:2008b,Baker:2005,Zhou:2006c}. Among them, GB models are heuristic approaches to polar solvation energy analysis. PC models resort to quantum mechanical calculations of induced solute charges. PB methods can be formally derived from Maxwell equations and statistical mechanics  for electrolyte solutions \cite{Beglov:1996,Netz:2000a,Holm:2001} and therefore offer the promise of  handling large biomolecules with sufficient accuracy and robustness \cite{David:2000,Onufriev:2000,Bashford:2000}.

Conceptually, the separation between continuum solvent and the discrete (atomistic) solute introduces an interface definition.   This  definition may take the form of analytic functions \cite{Zap,Grant:1995,Grant:2007} or nonsmooth boundaries  dividing the solute-solvent domains. The van der Waals surface, solvent accessible surface \cite{Lee:1971}, and molecular surface (MS) \cite{Richards:1977} are devised for this purpose and have found their success  in biophysical calculations \cite{Spolar:1989,Livingstone:1991,Crowley:2005,Kuhn:1992,Bergstrom:2003,Dragan:2004,Jackson:1995,Licata:1997}. It has been noticed that the performance of implicit solvent models is very sensitive to the interface definition  \cite{Dong:2003,Dong:2006,Nina:1999,Swanson:2005a}. This comes as no surprise because  many of these popular interface definitions are \textit{ad hoc} divisions of the solute and solvent domains based on rigid molecular geometry and neglecting  solute-solvent energetic interactions. Additionally, geometric singularities \cite{Connolly:1983,Sanner:1996}  associated with these surface definitions incur enormous computational instability \cite{Zhou:2006c,Yu:2007,Yu:2007a} and lead to conceptual difficulty in interpreting the sharp interface \cite{ZhanChen:2011a}.

The differential geometry (DG) theory of surfaces \cite{Willmore:1997} and associated geometric partial differential equations (PDEs)  provide a natural description of the solvent-solute interface. In 2005, Wei and his collaborators introduced curvature-controlled PDEs for generating molecular surfaces  in solvation analysis \cite{Wei:2005}.  The first variational solvent-solute interface, namely, the minimal molecular surface (MMS), was constructed in 2006 by Wei and coworkers  based on the DG theory of surfaces \cite{Bates:2006,Bates:2006f,Bates:2008}. MMSs are constructed by solving the  mean curvature flow, or the Laplace-Beltrami flow, and have been applied to the calculation of  electrostatic potentials and  solvation free energies  \cite{ZhanChen:2012, Bates:2008}. This approach was generalized to potential-driven geometric flows, which admit physical interactions, for the surface generation of  biomolecules in solution \cite{Bates:2009}. While our approaches were employed and/or modified by many others \cite{Cheng:2007e,Yu:2008g,SZhao:2011a,SZhao:2014a} for molecular surface and solvation analysis, our geometric PDE  \cite{Wei:2005} and variational surface models  \cite{Bates:2006,Bates:2008,Bates:2009} are, to our knowledge, the first of their kind for solvent-solute interface and solvation modeling.

Since the surface area minimization is  equivalent to the minimization of surface free energies, due to a constant surface tension, this approach can be easily incorporated into  the variational formulation of the PB theory \cite{Sharp:1990, Gilson:1993} to result in  DG-based full solvation  models  \cite{ZhanChen:2010a, Wei:1999}, following a similar approach by Dzubiella {\it et al} \cite{Dzubiella:2006,SGZhou:2015}. Our DG-based solvation models have been  implemented in the Eulerian formulation, where the solvent-solute interface is embedded in the three-dimensional (3D) Euclidean space and behaves like a smooth characteristic function  \cite{ZhanChen:2010a}. The resulting interface and associated dielectric function vary smoothly from their values in the solute domain to those in the solvent domain and are computationally robust.  An alternative implementation is the Lagrangian formulation \cite{ZhanChen:2010b} in which the solvent-solute boundary is extracted as a sharp surface at a given isovalue and subsequently used in the solvation analysis, including nonpolar and polar modeling.

One major advantage of our DG based solvation model is that it  enables the  synergistic coupling between the solute and solvent domains via the variation procedure. As a result,  our DG based solvation model is able to  significantly reduce the number of free parameters that users must ``fit'' or adjust in applications to real-world systems \cite{Thomas:2013}. 
It has been demonstrated  that physical parameters, i.e., pressure and surface tension obtained from experimental data, can be directly employed in our DG-based solvation models for accurate solvation energy prediction \cite{Daily:2013}. Another advantage  of our DG based solvation model is that it avoids the use of {\it ad hoc}   surface definitions   and its interfaces, particularly ones generated from the Eulerian formulation \cite{ZhanChen:2010a}, are free of troublesome geometric singularities that commonly occur in conventional solvent-accessible and solvent-excluded surfaces \cite{Connolly85,Sanner:1996}. As a result, our DG based solvation model bypasses the sophisticated interface techniques required for solving the PB equation  \cite{Yu:2007,Yu:2007a,Geng:2007a}. In particular, the smooth solvent-solute interface obtained from the Eulerian formulation \cite{ZhanChen:2010a} can be directly interpreted as the physical solvent-solute boundary profile. Additionally, the resulting smooth dielectric boundary can also have a straightforward physical interpretation. The other advantage of our DG based solvation model is that it is nature and easy to incorporate the density functional theory (DFT) in its variational formulation. Consequently,  it is able to reevaluate  and reassign the solute charge induced by solvent polarization effect during the solvation process \cite{ZhanChen:2011a}. The resulting  total energy minimization process  recreates or resembles the solvent-solute interactions, i.e., polarization,  dispersion, and polar and nonpolar coupling in a realistic solvation process.   Recently, DG based solvation model has been extended to DG based multiscale models for non-equilibrium processes in biomolecular systems \cite{Wei:1999,  Wei:2012, Wei:2013,DuanChen:2012a, DuanChen:2012b}. These models recover the DG based solvation model at the equilibrium  \cite{Wei:2012}.

Recently, we have demonstrated \cite{ZhanChen:2012} that the DG based nonpolar solvation model is able to  outperform many other methods   \cite{Gallicchio:2000,Wagoner:2006,Ratkova:2010}  in solvation energy predictions  for a large number nonpolar molecules.  The  root mean square error (RMSE) of our predictions was below 0.4kcal/mol,  which clearly indicates the potential power of the DG based solvation formulation. However, the  DG based full solvation model has not shown a similar superiority in accuracy, although it works very well \cite{ZhanChen:2010a,ZhanChen:2010b}. Having so many aforementioned advantages, our DG based solvation models ought to outperform other methods with a similar level of approximations. One obstacle that hinders the performance of our DG based {\it full} solvation model is the numerical instability in solving two strongly coupled and highly nonlinear PDEs, namely, the generalized Laplace-Beltrami  (GLB) equation and the generalized  PB (GPB) equation.  To avoid such instability, a strong parameter constraint was applied to the nonpolar part in our earlier work   \cite{ZhanChen:2010a,ZhanChen:2010b},  which results in the reduction of our model accuracy.

The objective of the present work is to explore a better parameter optimization of our DG based solvation models. A pair of conditions is prescribed to ensure the physical solution of the GLB equation, which leads to the well-posedness of the GPB equation. Such a well-posedness in turn renders the stability of solving the GLB equation. The stable solution of the coupled GLB and GPB equation enables us to optimize the model parameters and produce the highly accurate prediction of solvation free energies.   Some of the best results are obtained in the solvation free energy prediction of more than a hundred molecules of both polar and nonpolar types.

The rest of this paper is organized as the follows. To establish the notation and facilitate further development, we present a brief review of our DG based solvation models in Section \ref{Theory}.  By using the variational principle, we derive the coupled GLB and GPB equations. Necessary boundary conditions and initial values are prescribed to make this coupled system well-posed.  Section \ref{Algorithm} is devoted to parameter learning algorithms. We develop a protocol  to stabilize the iterative solution process of coupled nonlinear PDEs. We introduce perturbation and convex optimization methods to ensure stability of the numerical solution of the GLB equation in coupling with the GPB equation.    The newly achieved stability in solving the coupled PDEs leads to an appropriate minimization of  solvation free energies with respect to our model parameters. In Section \ref{Numerical-Result}, we show that for more than a hundred of compounds of various types, including both polar and nonpolar molecules, the present DG solvation model offers some of the most accurate solvation free energy prediction with the overall RMSE of 0.5kcal/mol.  This paper ends with a conclusion.

\section{The DG based solvation model}\label{Theory}

The free energy functional for our DG based  full solvation model can be expressed as  \cite{Wei:2009,ZhanChen:2010a,ZhanChen:2010b}
\begin{eqnarray} \label{eq8tot}
	\begin{aligned}
		G[S,\Phi]& =  \int\left\{ \gamma |\nabla S | +   p S   +   (1-S)U
		+S \left[  -\frac{\epsilon_m}{2}|\nabla\Phi|^2 + \Phi\ \rho_m\right] \right. \\
		& \left. +(1-S)\left[-\frac{\epsilon_s}{2}|\nabla\Phi|^2-k_B T \sum_{\alpha} \rho_{\alpha 0}\left( e^{-\frac{q_{\alpha }\Phi }{k_B T }}-1\right) \right] \right\} d{\bf{r}}, \quad {\bf r} \in {\mathbb R}^3
	\end{aligned}
\end{eqnarray}
where $\gamma$ is the surface tension, $p$ is the hydrodynamic pressure difference between solvent and solute, and $U$ denotes the solvent-solute non-electrostatic interactions represented by the Lennard-Jones potentials in the present work. Here  $0 \leq S \leq1$ is a hypersurface or simply surface function that characterizes the solute domain and embeds the 2D surface in ${\mathbb R}^3$, whereas $1-S$ characterizes the solvent domain \cite{ZhanChen:2010a}. Additionally, $\Phi$ is the electrostatic potential and $\epsilon_s$ and $\epsilon_m$ are the dielectric constants of the solvent and solute, respectively.  Here $k_B$ is the Boltzmann constant, $T$ is the temperature, $\rho_{\alpha 0}$ denotes the reference bulk concentration of the $\alpha$th solvent species, and $q_{\alpha }$ denotes the charge valence of the $\alpha$th solvent species, which is zero for an uncharged solvent component. We use $\rho_m$ to represent the charge density of the solute. The charge density is often modeled by a point charge approximation
$$\rho_m=\sum_{j}^{N_m}Q_{j} \delta ({\bf{r}}-{\bf{r}}_{j}),
$$
where $Q_{j}$ denoting the partial charge of the $j$th atom in the solute. Alternatively, the charge density computed from the DFT, which  changes during the iteration or energy minimization,  can be directly employed as well   \cite{ZhanChen:2011a}.

In Eq. (\ref{eq8tot}), the first three terms consist of  the so called nonpolar solvation free energy functional while the last two terms form  the polar one.
After the variation with respect to $S$,  we construct the following generalized Laplace-Beltrami (GLB) equation by using a procedure discussed in our earlier work \cite{Bates:2009}
\begin{eqnarray}\label{eq10surf}
	\frac{\partial S}{\partial t}&=& |\nabla S |\left[\nabla\cdot\left(\gamma\frac{\nabla S}{|\nabla S |}\right)
   + V\right],
\end{eqnarray}
where the potential driven term is given by
$$	V=  -p  +  U
	+\frac{\epsilon_m}{2}|\nabla\Phi|^2 - \Phi\ \rho_m-\frac{\epsilon_s}{2}|\nabla\Phi|^2 - k_B T \sum_{\alpha} \rho_{\alpha 0}\left( e^{-\frac{q_{\alpha }\Phi }{k_B T }}-1\right).
$$
As in the nonpolar case, solving the generalized Laplace-Beltrami equation (\ref{eq10surf}) generates the solvent-solute interface through the surface function $S$.

Additionally, variation with respect to $\Phi$ gives rise to the generalized Poisson-Boltzmann (GPB) equation:
\begin{eqnarray}\label{eq13poisson}
	-\nabla\cdot\left(\epsilon(S) \nabla\Phi\right)= S\rho_m
    +(1-S)\sum_{\alpha} q_{\alpha}\rho_{\alpha 0}e^{-\frac{q_{\alpha }\Phi }{k_B T }}, 
\end{eqnarray}
where $\epsilon(S)=(1-S)\epsilon_s+S\epsilon_m$ is the generalized permittivity function.
As shown in our earlier work \cite{Wei:2009,ZhanChen:2010a}, $\epsilon(S)$ is a smooth dielectric function gradually varying from $\epsilon_m$ to $\epsilon_s$. Thus, the solution procedure of the GPB equation avoids many numerical difficulties of solving elliptic equations with discontinuous coefficients \cite{Zhao:2004,Zhou:2006c,Zhou:2006d,Yu:2007c,Yu:2007a} in the standard PB equation.

The GLB (\ref{eq10surf}) and GBP (\ref{eq13poisson}) equations form a highly nonlinear system, in which the GLB equation is solved for the interface profile $S$ of the solute and solvent. The interface profile determines the dielectric function $\epsilon(S)$ in the GPB equation. The GPB equation is solved for the electrostatics potential $\Phi$ that behaves as an external potential in the GLB equation. The strongly coupled system should be solved in self-consistent  iterations.

For  GLB equation (\ref{eq10surf}), the computational domain is $\Omega/\Omega_m^{\rm vdW} $, where $\Omega_{m}^{\rm vdW}$ is the solute van der Waals domain given by $\Omega_m^{\rm vdW} = \bigcup_i B(r_i^{\rm vdW})$. Here $ B(r_i^{\rm vdW})$ is the $i$th ball in the solute centered at $\mathbf{r}_i$ with  van der Waals radius $r_i^{\rm vdW} $.
We apply the following Dirichlet  boundary condition to $S({\bf r},t)$
\begin{equation}
    \label{Dirichel-BC}
    S(\mathbf{r},t)=\left\{ \begin{array}{ll}
		0, &\forall \mathbf{r} \in \partial\Omega\\
		1, &\forall \mathbf{r} \in \partial\Omega_{m}^{\rm vdW}.
		\end{array} \right.
\end{equation}
The initial value of  $S({\bf r},t)$ is given by
\begin{equation}
    \label{S-initialV}
		S(\mathbf{r},0)=\left\{\begin{array}{ll}
		1, &\forall \mathbf{r} \in \partial\Omega_m^{\rm ext},\\
		0, & {\rm otherwise},
		\end{array}\right.
\end{equation}
where $\partial\Omega_m^{\rm ext}$ is the boundary of the extended solute domain constructed by $\Omega_m^{\rm ext} = \bigcup_i B(r_i^{\rm vdW}+r^{\rm probe})$. Here $ B(r_i^{\rm vdW}+r^{\rm probe})$  has an extended radius of $r_i^{\rm vdW}+r^{\rm probe}$ with  $r^{\rm probe}$ being the probe radius, which is set to 1.4\AA~ in the present work.

For  GPB equation (\ref{eq13poisson}),   the computational domain is $\Omega$. We  set the Dirichlet   boundary condition via the Debye-H\"uckel expression,
\begin{equation}
    \label{DH-BC}
    \Phi(\mathbf{r})=\sum_{i=1}^{N_m} \frac{Q_i}{\epsilon_s|\mathbf{r}-\mathbf{r}_i|} e^{-\bar{\kappa}|{{\bf r}-{\bf r}_i}|}, \quad \forall {\bf r}\in \partial\Omega,
    \end{equation}
where $\bar{\kappa}$ is the modified Debye-H\"uckel screening function \cite{ZhanChen:2010b}, which  is zero if  there is no salt molecule in the solvent.  Note that no interface condition \cite{Yu:2007} is needed as $S$ and $\epsilon(S)$ are smooth functions in general for $t>0$. Consequently, the resulting    GBP (\ref{eq13poisson}) equation is easy to solve.

To compare with experimental solvation data, one needs to compute the total solvation free energy, which, in our DG based solvation model,  is obtained as
\begin{equation}
\label{DGsolv-solvation}
\Delta G = \Delta G^{\rm P} + G^{\rm NP},
\end{equation}
where $\Delta G^{\rm P}$ is the    electrostatic solvation free energy,
\begin{equation}
\label{reaction-field-potential}
\Delta G^{\rm P}=\frac{1}{2}\sum_{i=1}^{N_m} Q_i\left[\Phi(\mathbf{r}_i)-\Phi_{h }(\mathbf{r}_i)\right]
\end{equation}
where   $\Phi_{h }$ is the solution of the above the GPB model in a homogenous system, obtained by setting a constant permittivity function $\epsilon(\mathbf{r})=\epsilon_m$ in the whole domain $\Omega$.
The  nonpolar energy $G^{\rm NP} $ is computed by
\begin{eqnarray} \label{eqnonpolar}
 G^{\rm NP}= \int \left[ \gamma |\nabla S | +   p S +   (1-S)U \right] d{\bf{r}}.
\end{eqnarray}

The DG based solvation model  is formulated as a coupled GLB and GPB equation system, in which the GLB equation provides the solvent solute boundary for solving the GPB, while the GPB equation produces the external potential in the GLB equation for the surface evolution. The solution procedure for this coupled system has been discussed in our earlier work   \cite{ZhanChen:2010a,ZhanChen:2010b}. Essentially, for the GLB equation, an alternating direction implicit (ADI) scheme is utilized for the time integral, in conjugation with   the second order finite difference method for the spatial discretization.  The GPB equation is discretized by a standard second order finite difference scheme and the resulting algebraic equation system is solved by using a standard Krylov subspace method based solver  \cite{ZhanChen:2010a,ZhanChen:2010b}.

\section{Parametrization methods and algorithms} \label{Algorithm}

 To solve the  above coupled equation system, a set of parameters that appeared in the GLB equation, namely, surface tension $\gamma$, hydrodynamic pressure difference $p$ and the product of solvent density $\rho_\alpha \varepsilon_j \doteq \tilde{\varepsilon}_{j\alpha}$, should be predetermined. Unfortunately, this coupled system   is unstable at the certain choices of parameters. Specifically, for certain $V$, one may have $S>1$ or $S<0$, which leads to unphysical $\epsilon(S)$ and unphysical solution of  GPB equation (\ref{eq13poisson}) and thus gives rise to a divergent $S$.  This instability can seriously reduce the model accuracy  \cite{ZhanChen:2010a,ZhanChen:2010b}.

For a concise description of our algorithm, we assume that there is only one solvent component (water) and denote the parameter set as:
\begin{equation}
\label{Parameters-set}
{\bf P}=\{\gamma, p, \tilde{\varepsilon}_{1}, \tilde{\varepsilon}_{2}, \cdots, \tilde{\varepsilon}_{N_T}\}
\end{equation}
where $N_T$ is the number of  types of atoms in the solute molecule.

As mentioned in the previous part, the parameter set  ${\bf P}$ used in solving the coupled PDEs should meet   two requirements, namely, the stability of solving the coupled PDEs and the optimal prediction of the solvation free energy (or fitting the experimental solvation free energy in the best approach). Based on these two criteria we introduce a two-stage numerical procedure to optimize the parameter set and solve the coupled PDEs:
\begin{itemize}
\item Explore the stability conditions of the coupled PDEs by introducing an auxiliary system via a small perturbation;

\item Optimize the parameter set by an iteratively scheme satisfying the stability constraint.
\end{itemize}

\subsection{Stability conditions}
In this part we investigate the stability conditions for the numerical solution to the coupled PDEs (\ref{eq10surf}) and (\ref{eq13poisson}). The basic idea is to utilize a small perturbation method. It is known that omitting  the external potential in the GLB  equation yields the Laplace-Beltrami (LB) equation:
\begin{equation}
\label{LBE}
\frac{\partial S}{\partial t}=|\nabla S|\nabla\cdot\left(\gamma\frac{\nabla S}{|\nabla S|}\right)
\end{equation}
This equation is of diffusion type and is well posed with the Dirichlet type of boundary conditions provided $\gamma>0$. Numerically it is easy to solve Eq. (\ref{LBE}) to  yield the profile of the solvent solute boundary.

After solving the LB equation   (\ref{LBE}), we use the generated smooth profile of the solvent solute boundary to determine the permittivity function in the GPB equation.  For simplicity,  we consider a pure water solvent,
\begin{equation}
\label{PBE-water}
-\nabla\cdot\left(\epsilon(S)\nabla\Phi\right)=S\rho_m.
\end{equation}
Without the external potential the system of Eqs. (\ref{LBE})-(\ref{PBE-water}) can be solved   stably by first solving the LB equation  and then the GPB equation.

Motivated by the above observation, if the external potential is dominated by the mean curvature term, the  stability of coupled GPB and GLB equations can be preserved.  Based on numerical experiments, the Lennard-Jones interaction between the solvent and solute is usually small since this term is constrained by the nonpolar free energy in our model. In our method, we enforce  the following constraint conditions to make the coupled system well-posed in the numerical sense
\begin{equation}
\label{gamma-constraint}
\gamma>\gamma_0>0,
\end{equation}
and
\begin{equation}
\label{pressure-constraint}
|p|\leq \beta \gamma,
\end{equation}
where $\gamma_0$ and $\beta$ are some appropriate positive constants.

In summary, the original problem is transformed into optimizing  parameters in the following system to attain the best solvation free energy fitting with experimental results:
\begin{equation}
\label{Constrained-PDE-optimization}
\left\{
  \begin{array}{ll}
    \frac{\partial S}{\partial t}=|\nabla S|\left[\nabla\cdot\left(\gamma\frac{\nabla S}{|\nabla S|} \right) -p+U+\frac{1}{2}\epsilon_m|\nabla\Phi|^2-\frac{1}{2}\epsilon_s|\nabla\Phi|^2 \right],  \\
    -\nabla\cdot\left(\epsilon(S)\nabla\Phi\right)=S\rho_m,  \\
    \gamma>\gamma_0>0,   \\
    |p|\leq \beta \gamma.
  \end{array}
\right.
\end{equation}
Note that the potential $\rho_m\Phi$ is omitted in the GLB equation (\ref{Constrained-PDE-optimization}),   because we have already enforced the Dirichlet boundary condition in the GLB equation, while $\rho_m$ is inside the van der Waals surface.

\begin{remark}
Based on large amount of numerical tests, it is found that there is no need to enforce the constraint conditions on the parameters that appear in the Lennard-Jones term.  When this term is used to fit the solvation energy with experimental results, the parameters can be bounded in a small neighborhood of 0 automatically during the fitting procedure.  These parameters essentially do not affect the numerical stability.
\end{remark}

\subsection{Self-consistent approach for solving the coupled PDEs}
In this part, we propose a self-consistent approach to solve the coupled GLB and GPB equations for a given set of parameters. Basically, the coupled system is solved iteratively until both the electrostatic solvation free energy $\Delta G^{\rm P}$ given in Eq. (\ref{reaction-field-potential}) and the surface function $S$ are both converged.  Here the surface function is said to be converged provided that the surface area and enclosed volume are both converged.

We present an algorithm for solving the following coupled systems:
\begin{equation}
\label{GPB2}
-\nabla\cdot(\epsilon(S)\nabla\Phi)=S\rho_m,
\end{equation}
and
\begin{equation}
\label{GLB2}
\frac{\partial S}{\partial t}=|\nabla S|\left[\nabla\cdot\left(\gamma\frac{\nabla S}{|\nabla S|}\right)+V_e\right],
\end{equation}
where $V_e$ is the external potential which is defined as:
\begin{itemize}
\item \textbf{Auxiliary system:} $V_e=\frac{1}{2}(\epsilon_m-\epsilon_s)|\nabla\Phi|^2$,
\item \textbf{Full system:} $V_e=-p+U+\frac{1}{2}(\epsilon_m-\epsilon_s)|\nabla\Phi|^2$.
\end{itemize}

Dirichlet boundary conditions are employed for both   GPB (\ref{GPB2}) and    GLB (\ref{GLB2}) equations with auxiliary and full external potentials, giving rise to  a well-posed coupled system. The smooth profile of the solvent-solute boundary enables the direct use of the second order central finite difference scheme to achieve the second order convergence in discretizing the GPB equation. The biconjugate gradient scheme is used to solve  the resulting algebraic equation system. The GLB equation of both the auxiliary and full systems can be solved by the central finite difference discretization of the spatial domain and the forward Euler time integrator  for the time domain discretization. 

\begin{remark}
For the sake of simplicity, in the current work, we employed the central finite difference scheme for spatial domain discretization in both GPB and GLB equations, and forward Euler integrator for the time domain discretization of GLB equation. For stability consideration, in the discretization of the GLB equation, the discretization step size of temporal and spatial domain  satisfies the Courant-Friedrichs-Lewy condition. To accelerate the numerical integration, a multigrid solver can be employed for GBP equation, and an alternating direction implicit scheme \cite{ZhanChen:2010a}, which is unconditionally stable, can be utilized for the temporal integration. However, detail discussion of these accelerated schemes is beyond the scope of the present work.
\end{remark}

A pseudo code is given in Algorithm  \ref{SCF-PB-LB}  to offer a general framework for solving the coupled GLB and GPB equations in a self-consistent manner. The  outer iteration controls the convergence of the GPB equation through measuring the change of  electrostatic solvation free energy in two  adjacent iterations, while the  inner iteration controls the convergence of the GLB equation based on the variation of surface areas and  enclosed volumes through  the surface function $S$.
The variables ${\rm \Delta G^{\rm P}_1}$, ${\rm \Delta G^{\rm P}_2}$, ${\rm Area_1}$, ${\rm Area_2}$, ${\rm Vol_1}$, and ${\rm Vol_2}$ denote the electrostatic solvation free energy, surface area, and volume enclosed by the surface of two immediate iterations, respectively.
\begin{algorithm}
\caption{Self-consistent algorithm for the coupled  GPB and GLB system}\label{SCF-PB-LB}
\begin{algorithmic}[1]
\Procedure{GPB-GLB-Solver}{}
\State \textbf{Initialize: ${\rm \Delta G^{\rm P}_1}=0$, ${\rm \Delta G^{\rm P}_2}=100$, ${\rm Area_1}=0$, ${\rm Area_2}=100$, ${\rm Vol_1}=0$, ${\rm Vol_2}=100$}
\State \textbf{do while} ($|{\rm \Delta G^{\rm P}_1}-{\rm \Delta G^{\rm P}_2}|<\epsilon_1$) 
\State \hskip 1.0cm ${\rm \Delta G^{\rm P}_1}\leftarrow{\rm \Delta G^{\rm P}_2}$
\State \hskip 1.0cm \textbf{do while} ($|{\rm Area_1}-{\rm Area_2}|<\epsilon_2$ .and. $|{\rm Vol_1}-{\rm Vol_2}|<\epsilon_3$)
\State \hskip 2.0cm ${\rm Area_1}\leftarrow{\rm Area_2}$, ${\rm Vol_1}\leftarrow{\rm Vol_2}$.
\State \hskip 2.0cm Update the surface profile function $S$ by solving the GLB equation  (\ref{GLB2}).
\State \hskip 2.0cm ${\rm Area_2}=\int_\Omega Sd\mathbf{r}$, ${\rm Vol_2}=\int_\Omega |\nabla S|d\mathbf{r}$.
\State \hskip 1.0cm \textbf{enddo}
\State \hskip 1.0cm Solve the GPB equation  (\ref{GPB2}) in both vacuum and solvent with the previous updated surface profile.
\State \hskip 1.0cm Update the polar solvation free energy $\Delta G^{\rm P}_2$ according to Eq. (\ref{reaction-field-potential}).
\State \textbf{enddo}
\EndProcedure
\end{algorithmic}
\end{algorithm}
The parameters $\epsilon_1, \epsilon_2$ and $\epsilon_3$ are the threshold constants and all set to $0.01$ in the current implementation.

\begin{remark}
In solving the GLB equation, during each updating, to ensure the stability, instead of the fully update, we update it partially, i.e., the updated solution is the weighted sum of the new solution of the current GLB solution $S_{\rm new}$ and the old solution of the GLB equation in the previous step $S_{\rm old}$:
\begin{equation}
\label{S-update}
S=a_1S_{\rm new}+(1-a_1)S_{\rm old},
\end{equation}
where $a_1$ is a constant and set to 0.5  in the present work.
\end{remark}

\subsection{Convex optimization for parameter learning}

In this part,  we present the parameter optimization scheme.  In our approach, parameters start from an initial guess and then are updated sequentially until reaching the convergence. Here the convergence is measured by the root mean square (RMS) error between the fitted and experimental solvation free energies for a given set of molecules.

Consider the parameter optimization for a given group of molecules, denoted as $\{T_1, T_2, \cdots, T_n\}$.  As discussed above the parameter set is ${\bf P}$. To optimize the parameter set ${\bf P}$, we start from   GPB equation (\ref{GPB2}) and the auxiliary system of   GLB equation  (\ref{GLB2}) with $\gamma=0.05$. After solving the initial coupled system by using Algorithm  \ref{SCF-PB-LB}, we obtain the following quantities for each molecule in the training set:
\begin{align}
\label{quantities-total-j}
\left\{{\rm \Delta G^P_j}, {\rm Area_j}, {\rm Vol_j},
 \left(\sum_{i=1}^{N_m} \delta_i^1 \int_{\Omega_s} \left[\left(\frac{\sigma_s+\sigma_1}{||\mathbf{r}-\mathbf{r}_i||}\right)^{12}-
2\left(\frac{\sigma_s+\sigma_1}{||\mathbf{r}-\mathbf{r}_i||}\right)^{6}\right] d\mathbf{r}\right)_j,\right.\\
 \left.
\cdots,\right.\\
 \left.
\left(\sum_{i=1}^{N_m} \delta_i^{N_T} \int_{\Omega_s} \left[\left(\frac{\sigma_{s}+\sigma_{N_T}}{||\mathbf{r}-\mathbf{r}_i||}\right)^{12}-
2\left(\frac{\sigma_{s}+\sigma_{N_T}}{||\mathbf{r}-\mathbf{r}_i||}\right)^{6}\right] d\mathbf{r}\right)_j \right\}
\end{align}
where $j=1, 2, \cdots,  n$. Here $N_m$ and $N_T$ denote the number of atoms  and types of atoms in a specific molecule. The last  few terms involve semi-discrete and semi-continuum Lennard-Jones potentials \cite{ZhanChen:2010a}.
Additionally,
$$
\delta_i^j=\left\{
             \begin{array}{ll}
               1, & \hbox{if atom $i$ belongs to type $j$,} \\
               0, & \hbox{otherwise.}
             \end{array}
           \right.
$$
where $i=1, 2, \cdots, N_m$; $j=1, 2, \cdots N_T$; $\sigma_i, i=1, 2, \cdots, N_T$ is the atomic radius of the $i$th type of atoms. Therefore, atoms of the same type have a common atomic radius and fitting parameter $\tilde{\varepsilon}$.

The predicted solvation free energy for molecule $j$ can be represented as:
\begin{eqnarray}
\label{pred-eng-j}
{\rm \Delta G_j}={\rm \Delta G^P_j}+\gamma {\rm Area_j}+p{\rm Vol_j}+\tilde{\varepsilon}_{1} \left(\sum_{i=1}^{N_m} \delta_i^1 \int_{\Omega_s} \left[\left(\frac{\sigma_s+\sigma_1}{||\mathbf{r}-\mathbf{r}_i||}\right)^{12}-
2\left(\frac{\sigma_s+\sigma_1}{||\mathbf{r}-\mathbf{r}_i||}\right)^{6}\right] d\mathbf{r}\right)_j \\
+\cdots+\tilde{\varepsilon}_{N_T} \left(\sum_{i=1}^{N_m} \delta_i^{N_T} \int_{\Omega_s} \left[\left(\frac{\sigma_{s}+\sigma_{N_T}}{||\mathbf{r}-\mathbf{r}_i||}\right)^{12}-
2\left(\frac{\sigma_s+\sigma_{N_T}}{||\mathbf{r}-\mathbf{r}_i||}\right)^{6}\right] d\mathbf{r}\right)_j.
\end{eqnarray}
We denote the predicted solvation free energy for the given set of molecules as $\Delta \mathbf{G}({\bf P})\doteq \left\{\Delta G_1, \Delta G_2, \cdots, \Delta G_n\right\}$, which is a function of the parameter  set ${\bf P}$, and denote the corresponding experimental solvation free energy as $\Delta \mathbf{G}^{\rm Exp}\doteq \left\{\Delta G^{\rm Exp1}, \Delta G^{\rm Exp2}, \cdots, \Delta G^{{\rm Exp}n} \right\}$.

Then the parameter optimization problem in the coupled PDEs given by Eqs. (\ref{Constrained-PDE-optimization}) can be transformed into the following regularized and constrained optimization problem:
\begin{equation}
\label{convex-minimization}
\min_{\bf P} \left( ||\Delta \mathbf{G}({\bf P})-\Delta \mathbf{G}^{\rm Exp}||_2+\lambda||{\bf P}||_2\right),
\end{equation}
s.t.
\begin{equation}
\label{constraint1}
\gamma\geq \gamma_0,
\end{equation}
and
\begin{equation}
\label{constraint2}
|p|\leq \beta \gamma,
\end{equation}
where $||*||_2$ is the $L_2$ norm of the quantity $*$ and  $\lambda$ is the regularization parameter chosen to be 10 in the present work to ensure the dominance of the first term and avoid overfitting.   Here $\gamma_0$ and $\beta$ are set respectively  to  $0.05$ and $0.1$ in the present implementation, which guarantees the stability of the coupled system according to a large amount of numerical tests.

It is obvious that the objective function  (\ref{convex-minimization}) in the optimization is a convex function, meanwhile the solution domain restricted by  constraints  (\ref{constraint1})-(\ref{constraint2}) forms a convex domain. Therefore the optimization problem given by Eqs. (\ref{convex-minimization})-(\ref{constraint2}) is a convex optimization problem, which was studied by Grant and Boyd \cite{cvx,  gb08}.

After solving the above convex optimization problem, parameter set ${\bf P}$ is updated and used again in solving the coupled GLB and GPB system, i.e., Eqs.  (\ref{GLB2}) and (\ref{GPB2}). Repeating the above procedure, a new group of predicted solvation free energies together with a new group of parameters is obtained. This procedure is repeated until the RMS error between the predicted and experimental solvation free energies in  two sequential iterations is within a given threshold.

\begin{figure}[!ht]
\small
\centering
\includegraphics[width=8cm,height=6cm]{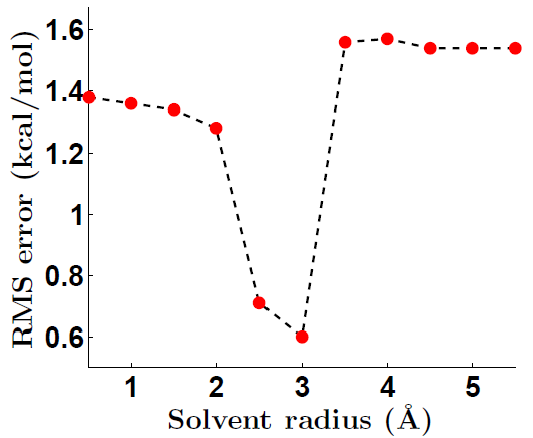}
\caption{The relation between the solvent radii and the RMS error of the SAMPL0 test set. The local minimum appears at the solvent radii 3.0 \AA, with RMS error being 0.6 kcal/mol for a set of 17 molecules.}
\label{solvent-radii-RMS}
\end{figure}

\subsection{Algorithm for parameter optimization and solution of the coupled PDEs}
Based on the preparation made in the previous two subsections, namely, the self-consistent approach for solving the coupled GLB and GPB system and the parameter optimization, we  provide the combined algorithm for the parameter  optimization and solving the coupled system  for a given set of molecules.

Algorithm \ref{parameters-learning} offers a parameter learning pseudo code for a given group of molecules. This algorithm is formulated by  combining outer and inner self-consistent iterations.  The outer iteration controls the convergence of the optimized parameters via two controlling parameters, ${\rm Err}_1$ and ${\rm Err}_2$, denoting the RMS error  between predicted and experimental solvation free energies in   two sequential iterations. The inner iteration implements the solution to the GLB and GPB equations by Algorithm \ref{SCF-PB-LB}.

\begin{algorithm}
\caption{Parameters learning for a given group of molecules}\label{parameters-learning}
\begin{algorithmic}[1]
\Procedure{Parameters-Learning}{}
\State \textbf{Initialize}: ${\rm Err}_1=0$, ${\rm Err}_2=100$
\State Solve the coupled GPB and GLB system, where GLB utilizes the auxiliary equation (\ref{GLB2}).
\State Solve the constrained optimization problem Eqs. (\ref{convex-minimization})-(\ref{constraint2}) to obtain the initial parameter  set ${\bf P}_0$.
\State Update ${\rm Err}_1$ to be the RMS error between experimental and predict results in the above step.

\State \textbf{do while} ($|{\rm Err}_1-{\rm Err}_2|<\epsilon_4$) 
\State \hskip 1.0cm ${\rm Err}_2 \leftarrow {\rm Err}_1$.
\State \hskip 1.0cm Solve the coupled GPB and GLB system, where GLB system with parameters set ${\bf P}_0$.
\State \hskip 1.0cm Solve the constrained optimization problem Eqs. (\ref{convex-minimization})-(\ref{constraint2}) to get the updated parameters set ${\bf P}$.
\State \hskip 1.0cm Update ${\rm Err}_1$ to be RMS error between experimental and predict results in the previous optimization step.
\State \hskip 1.0cm Update ${\bf P}_0\leftarrow {\bf P}$.
\State \textbf{enddo}
\EndProcedure
\end{algorithmic}
\end{algorithm}
The threshold parameter $\epsilon_4$ is set to   $0.01$ in the present work.

\section{Numerical results} \label{Numerical-Result}

In this section we present the numerical study of the DG based solvation model using the proposed parameter optimization algorithms. We first explore the optimal solvent radius used in  the van der Waals interactions. Due to the high nonlinearity, the solvent radius cannot be automatically optimized and its optimal value is obtained via searching the parameter domain. We show that for a group of molecules, there is a local  minimum in the RMS error when the solvent radius is varied. The corresponding optimal solvent radius  is adopted for other molecules. Additionally, we consider a large number of molecules with known experimental solvation free energies to test  the proposed parameter optimization algorithms. These molecules are of both polar and nonpolar types and are divided into six groups: the SAMPL0 test set \cite{Nicholls:2008}, the alkane, alkene, ether, alcohol and phenol types \cite{Mobley:2014}. It is found that our DG based solvation model works really well for these molecules.   Finally, to demonstrate the predictive power of the present DG based solvation model, we perform a five-fold cross validation \cite{EMSL:2009} for alkane, alkene, ether, alcohol and phenol types of molecules. It is found that  training and validation errors are of  the same level, which confirms the ability of our model for the solvation free energy prediction.

The SAMPL0 molecule structural conformations are adopted from the literature  with ZAP 9 radii and  the OpenEye-AM1-BCC v1 charges   \cite{Nicholls:2008}. For other molecules, structural conformations are obtained from FreeSolv \cite{Mobley:2014}. Amber GAFF force field is utilized for the charge assignment \cite{AMBER15}. The van der Waals  radii as well as the atomic radii of Hydrogen, Carbon and Oxygen atoms are set to  1.2, 1.7 and 1.5\AA, respectively. The grid spacing is set to 0.25\AA~ in all of our calculations (discretization and integration). The computational domain is set to  the bounding box of the solute molecule with an extra buffer length of 6.0 \AA.

\begin{figure}[!ht]
\small
\centering
\includegraphics[width=8cm,height=6cm]{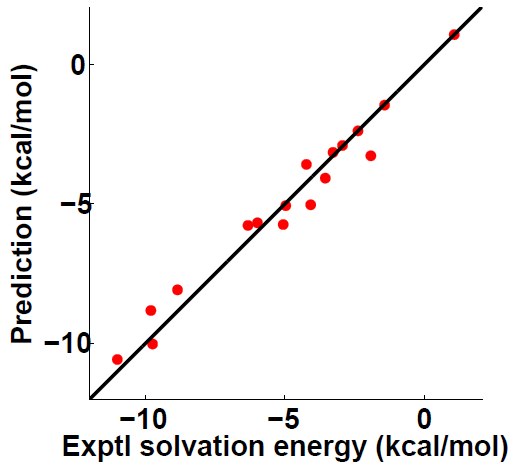}
\caption{The predicted and experimental solvation free energy for the 17 molecules in the SAMPL0 test set.}
\label{SAMPL0}
\end{figure}

\begin{table}[!ht]
\centering
\caption{The solvation free energy prediction for the SAMPL0 set.  Energy is in the unit of kcal/mol.}
\begin{tabular}{rrrrrr}
\cline{1-6}
Name  &$\Delta G^{\rm P}$ &$ G^{\rm NP}$ &${\Delta G}$ &$\Delta G^{\rm Exp}$\cite{Nicholls:2008} & Error \\
\hline
Glycerol triacetate              &-10.60    &2.53      &-8.07      &-8.84      &-0.77     \\
Benzyl bromide                   &-4.31     &1.93      &-2.38      &-2.38      &0.00     \\
Benzyl chloride                  &-4.45     &1.18      &-3.27      &-1.93      &1.34     \\
m-Bis (trifluoromethyl) benzene  &-2.62     &3.70      &1.08       &1.07       &-0.01     \\
N,N-Dimethyl-p-methoxybenzamide  &-8.35     &-2.22     &-10.57     &-11.01     &-0.45     \\
N,N-4-Trimethylbenzamide         &-6.93     &-3.09     &-10.03     &-9.76      &0.27     \\
bis-2-Chloroethyl ether          &-3.73     &-0.14     &-3.59      &-4.23      &-0.64     \\
1,1-Diacetoxyethane              &-7.07     &2.00      &-5.07      &-4.97      &0.10     \\
1,1-Diethoxyethane               &-3.58     &0.43      &-3.15      &-3.28      &-0.13     \\
1,4-Dioxane                      &-5.36     &-0.38      &-5.74      &-5.05      &0.69     \\
Diethyl propanedioate            &-7.07     &1.40      &-5.67      &-6.00      &-0.33     \\
Dimethoxymethane                 &-4.09     &1.19      &-2.90      &-2.93      &-0.03     \\
Ethylene glycol diacetate        &-7.66     &1.90      &-5.76      &-6.34      &-0.58     \\
1,2-Diethoxyethane               &-3.64     &0.45      &-4.09      &-3.54      &0.55     \\
Diethyl sulfide                  &-2.21     &0.76      &-1.47      &-1.43      &0.04     \\
Phenyl formate                   &-7.10     &2.08      &-5.02      &-4.08      &0.94     \\
Imidazole                        &-11.54    &2.71      &-8.83      &-9.81      &-0.98     \\
\hline
RMS                              &    &      &      &      &0.60   \\
\hline
\end{tabular}
\label{SAMPL0-table}
\end{table}

\subsection{Solvent radius}

In the present method,  the van der Waals radii of solute atoms are employed to define the van der Waals surface, which is used for setting up the boundary condition for the GLB equation. Additionally, solvent and solute atomic radii are used in the Lennard-Jones potentials. Atomic radii of solute are set to the van der Waals radii in the present work, whereas the solvent radius is considered an optimization parameter. We utilize a brute force approach for the solvent radii selection. We employ the SAMPL0 test set \cite{Nicholls:2008} as a benchmark. The solvent radius is varied from 0.5 \AA~ to 5.5 \AA~ away from van der Waals surface.  Due to the fast decay property of the Lennard-Jones interactions, the above setting enables the full inclusion of the Lennard-Jones interactions in our model.

Figure \ref{solvent-radii-RMS} depicts the RMS error of the 17 molecules from the SAMPL0 set at the different solvent radii calculated from the present DG based  solvation model. The result clearly demonstrates that  with the increase of the solvent radius, the RMS error decreases dramatically initially. The minimum appears at 3.0 \AA.  The further increase of the solvent radius leads to a rapid jump in the RMS error before it stabilizes around 1.54 kcal/mol. It is noted that  3.0 \AA~ is much larger than the commonly used solvent radius of 1.4 \AA~ in   Poisson-Boltzmann based implicit solvent models.  However, unlike in the commonly used implicit solvent models, the Lennard-Jones potential in our DG based solvation model is of half discrete and half continuum. Therefore, the solvent radius is an on-grid average value in the DG based solvation model.  In all the following computations, the solvent radius is set to 3.0 \AA.

\begin{figure}[!ht]
\small
\centering
\includegraphics[width=8cm,height=6cm]{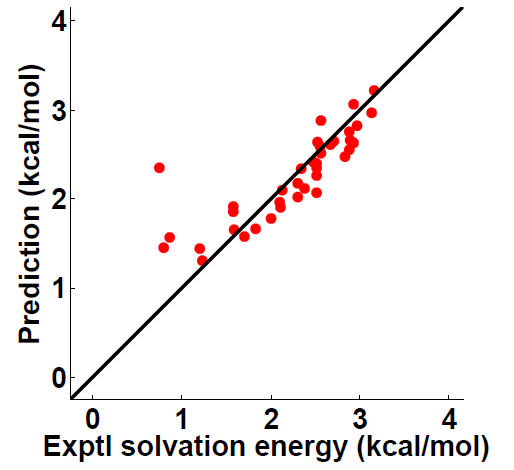}
\caption{The predicted and experimental solvation free energies for 38 alkane molecules.}
\label{pred-exp-alkane}
\end{figure}

\begin{figure}[!ht]
\small
\centering
\includegraphics[width=8cm,height=6cm]{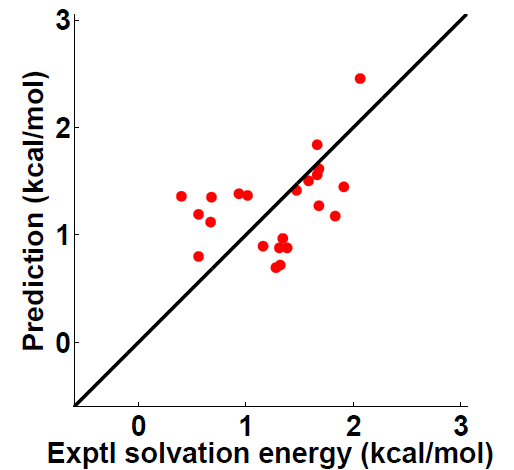}
\caption{The predicted and experimental solvation free energies for  22 alkene molecules.}
\label{pred-exp-alkene}
\end{figure}

\begin{table}[!ht]
\centering
\caption{The solvation free energy prediction for the alkane set.  All energies are in the unit of kcal/mol.}
\begin{tabular}{rrrrrr}
\cline{1-6}
Name  &$\Delta G^{\rm P}$ &$ G^{\rm NP}$ &${\Delta G}$ &$\Delta G^{\rm Exp}$\cite{Mobley:2014} & Error \\
\hline
octane                       &-0.13    &2.89      &2.76     &2.88      &0.12     \\
ethane                       &-0.04    &1.70      &1.66     &1.83      &0.17     \\
propane                      &-0.05    &1.83      &1.78     &2.00      &0.22     \\
cyclopropane                 &-0.08    &2.43      &2.35     &0.75      &-1.60     \\
isobutane                    &-0.07    &2.09      &2.02     &2.30      &0.28     \\
2,2-dimethylbutane           &-0.07    &2.34      &2.27     &2.51      &0.24     \\
isopentane                   &-0.07    &2.19      &2.12     &2.38      &0.26     \\
2,3-dimethylbutane           &-0.07    &2.41      &2.34     &2.34      &0.00     \\
3-methylpentane              &-0.08    &2.43      &2.35     &2.51      &0.16     \\
methylcyclopentane           &-0.10    &1.76      &1.66     &1.59      &-0.07     \\
n-butane                     &-0.07    &2.03      &1.96     &2.10      &0.14     \\
isohexane                    &-0.09    &2.49      &2.40     &2.51      &0.11     \\
2,4-dimethylpentane          &-0.09    &2.57      &2.48     &2.83      &0.35     \\
methylcyclohexane            &-0.10    &1.68      &1.58     &1.70      &0.12     \\
n-pentane                    &-0.08    &2.25      &2.17     &2.30      &0.13     \\
hexane                       &-0.09    &2.51      &2.42     &2.48      &0.06     \\
cyclohexane                  &-0.10    &1.40      &1.30     &1.23      &-0.07     \\
nonane                       &-0.14    &3.11      &2.97     &3.13      &0.16     \\
heptane                      &-0.11    &2.73      &2.62     &2.67      &0.05     \\
cyclopentane                 &-0.10    &1.54      &1.44     &1.20      &-0.24     \\
cycloheptane                 &-0.11    &1.56      &1.45     &0.80      &-0.65     \\
cyclooctane                  &-0.12    &1.69      &1.57     &0.86      &-0.71     \\
neopentane                   &-0.06    &2.13      &2.07     &2.51      &0.44     \\
2,2,4-trimethylpentane       &-0.08    &2.74      &2.66     &2.89      &0.23     \\
3,3-dimethylpentane          &-0.07    &2.58      &2.51     &2.56      &0.05     \\
2,3-dimethylpentane          &-0.08    &2.72      &2.64     &2.52      &-0.12     \\
2,3,4-trimethylpentane       &-0.08    &2.96      &2.88     &2.56      &-0.32     \\
1,2-dimethylcyclohexane      &-0.10    &2.02      &1.92     &1.58      &-0.34     \\
3-methylhexane               &-0.09    &2.74      &2.65     &2.71      &0.06     \\
3-methylheptane              &-0.11    &2.94      &2.83     &2.97      &0.14     \\
1,4-dimethylcyclohexane      &-0.11    &2.02      &1.91     &2.11      &0.20     \\
2,2-dimethylpentane          &-0.08    &2.64      &2.56     &2.88      &0.32     \\
2-methylhexane               &-0.10    &2.73      &2.63     &2.93      &0.30     \\
decane                       &-0.16    &3.37      &3.21     &3.16      &-0.06     \\
propylcyclopentane           &-0.12    &2.21      &2.09     &2.13      &0.03     \\
cis-1,2-Dimethylcyclohexane  &-0.09    &1.95      &1.86     &1.58      &-0.28     \\
2,2,5-trimethylhexane        &-0.09    &3.15      &3.06     &2.93      &-0.13     \\
pentylcyclopentane           &-0.15    &2.73      &2.58     &2.55      &-0.04     \\
\hline
RMS                          &    &      &      &      &0.36   \\
\hline
\end{tabular}
\label{alkanes}
\end{table}

\subsection{Optimization results}

In this section, we illustrate the performance of our parameter optimization algorithms.  First, we provide the regression results of the SAMPL0 test set \cite{Nicholls:2008}.  Figure \ref{SAMPL0} shows the predicted and experimental solvation free energies based on the present model and optimization method. It is obvious that    predicted solvation free energies are highly consistent with the experimental ones. The  RMS error is 0.60 kcal/mol.

Table \ref{SAMPL0-table} shows the breakup of polar, nonpolar and total predicted solvation free energies. The experimental values and errors are also provided \cite{Nicholls:2008}.

Compared to our earlier prediction \cite{ZhanChen:2010a} in which the same model is employed but the parameters were not optimized in the present manner, the RMS error decreases dramatically from previous 1.76 kcal/mol to 0.60 kcal/mol for the same test set.  Note that the present RMS error (0.60 kcal/mol) is also significantly smaller than that of the explicit solvent approach (1.71 $\pm$ 0.05 kcal/mol) and that obtained by the PB based prediction  (1.87 kcal/mol) under the same structure, charge and radius setting \cite{Nicholls:2008}.  The present results confirm  the efficiency of the proposed new  parameter optimization algorithms    and demonstrate the accuracy and power of our DG based solvation models.

\begin{table}[!ht]
\centering
\caption{The solvation free energy prediction for the alkene set. All energies are in the unit of kcal/mol.}
\begin{tabular}{rrrrrr}
\cline{1-6}
Name  &$\Delta G^{\rm P}$ &$ G^{\rm NP}$ &${\Delta G}$ &$\Delta G^{\rm Exp}$\cite{Mobley:2014} & Error \\
\hline
ethylene                    &-0.27    &0.96      &0.69     &1.28      &0.59     \\
isoprene                    &-0.62    &1.97      &1.35     &0.68      &-0.67     \\
but-1-ene                   &-0.29    &1.17      &0.88     &1.38      &0.50     \\
butadiene                   &-0.56    &1.75      &1.19     &0.56      &-0.63     \\
pent-1-ene                  &-0.30    &1.57      &1.27     &1.68      &0.41     \\
prop-1-ene                  &-0.32    &1.03      &0.71     &1.32      &0.61     \\
2-methylprop-1-ene          &-0.37    &1.26      &0.89     &1.16      &0.27     \\
cyclopentene                &-0.37    &1.17      &0.79     &0.56      &-0.23     \\
2-methylbut-2-ene           &-0.40    &1.28      &0.87     &1.31      &0.44     \\
2,3-dimethylbuta-1,3-diene  &-0.65    &2.01      &1.36     &0.40      &-0.95     \\
3-methylbut-1-ene           &-0.27    &1.45      &1.18     &1.83      &0.65     \\
1-methylcyclohexene         &-0.38    &1.50      &1.11     &0.67      &-0.45     \\
penta-1,4-diene             &-0.53    &1.91      &1.38     &0.93      &-0.45     \\
hex-1-ene                   &-0.30    &1.81      &1.50     &1.58      &0.08    \\
hexa-1,5-diene              &-0.51    &1.88      &1.37     &1.01      &-0.36     \\
hept-1-ene                  &-0.33    &2.17      &1.84     &1.66      &-0.18     \\
hept-2-ene                  &-0.34    &1.96      &1.62     &1.68      &0.06     \\
4-Methyl-1-pentene          &-0.26    &1.71      &1.45     &1.91      &0.46     \\
2-methylpent-1-ene          &-0.33    &1.75      &1.42     &1.47      &0.05     \\
non-1-ene                   &-0.36    &2.81      &2.45     &2.06      &-0.39     \\
trans-2-Heptene             &-0.34    &1.90      &1.56     &1.66      &0.10     \\
trans-2-Pentene             &-0.30    &1.26      &0.96     &1.34      &0.38     \\
\hline
RMS                         &    &      &      &      &0.46   \\
\hline
\end{tabular}
\label{alkenes}
\end{table}

Additionally, we investigate the solvation free energies prediction of two families of nonpolar molecules, alkane and alkene, which were studied previous by using our DG based nonpolar solvation model \cite{ZhanChen:2012}. In the following, we demonstrate that the present DG based full solvation model can provide the same level of accuracy in the solvation free energy prediction for alkane and alkene molecules.

Figures  \ref{pred-exp-alkane}  and  \ref{pred-exp-alkene}  depict the predicted and experimental solvation free energies for 38 alkane and 22 alkene molecules, respectively. Tables \ref{alkanes}  and  \ref{alkenes}  list the polar, nonpolar, total and experimental solvation free energies for both families of solute molecules, respectively.  Except for one alkane molecule, namely, cycloprotane, whose predicted error  is 1.60 kcal/mol, the errors for all other  molecules are within 1 kcal/mol. The RMS errors of these two families are 0.36 and 0.46 kcal/mol, respectively.
 This level of accuracy is similar to our earlier results obtained by using our DG based nonpolar  solvation model \cite{ZhanChen:2012}, which does not involve the electrostatic (polar) model and is computationally easier to optimize.

\begin{figure}[!ht]
\small
\centering
\includegraphics[width=8cm,height=6cm]{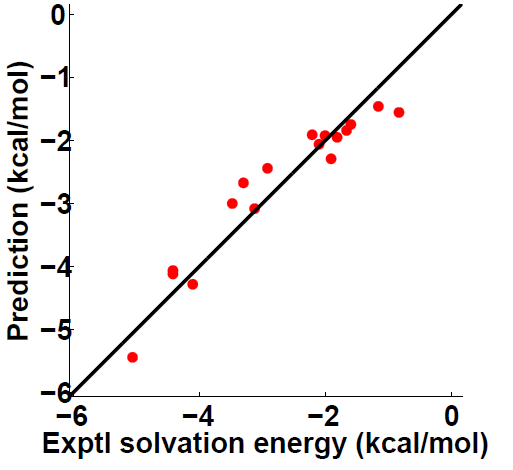}
\caption{The predicted and experimental solvation free energy for the 17 ether molecules.}
\label{pred-exp-ether}
\end{figure}

\begin{figure}[!ht]
\small
\centering
\includegraphics[width=8cm,height=6cm]{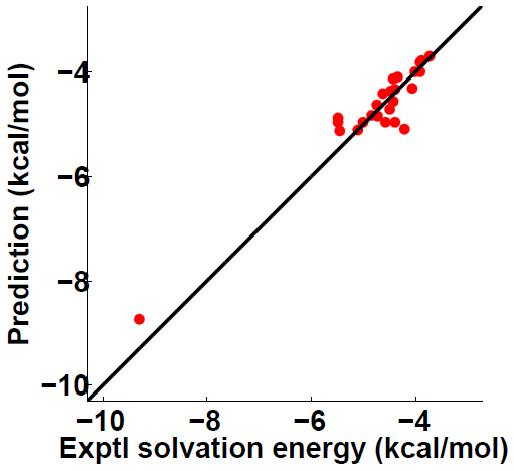}
\caption{The predicted and experimental solvation free energy for the 25 alcohol molecules.}
\label{pred-exp-alcohol}
\end{figure}

\begin{figure}[!ht]
\small
\centering
\includegraphics[width=8cm,height=6cm]{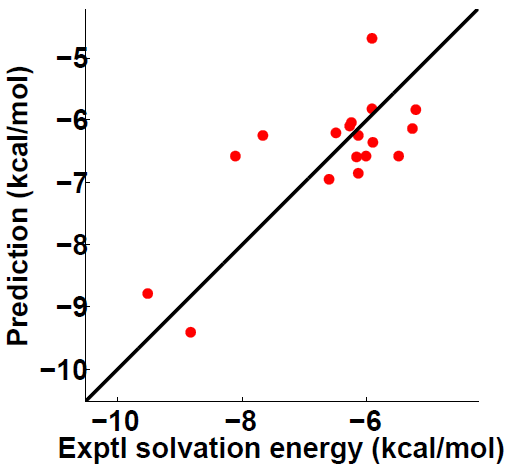}
\caption{The predicted and experimental solvation free energy for the 18 phenol molecules.}
\label{pred-exp-phenol}
\end{figure}

\begin{table}[!ht]
\centering
\caption{The solvation free energy prediction for the ether set. All energies are in the unit of kcal/mol.}
\begin{tabular}{rrrrrr}
\cline{1-6}
Name  &$\Delta G^{\rm P}$ &$ G^{\rm NP}$ &${\Delta G}$ &$\Delta G^{\rm Exp}$\cite{Mobley:2014} & Error \\
\hline
ethoxyethane                    &-4.08    &2.33      &-1.75     &-1.59      &0.16     \\
2-methyltetrahydrofuran         &-4.10    &1.43      &-2.67     &-3.30      &-0.63     \\
tetrahydrofuran                 &-4.36    &1.36      &-3.00     &-3.47      &-0.47     \\
1-propoxypropane                &-3.75    &2.29      &-1.46     &-1.16      &0.30     \\
methoxymethane                  &-4.55    &2.26      &-2.29     &-1.91      &0.36     \\
tetrahydropyran                 &-4.17    &1.09      &-3.07     &-3.12      &-0.05     \\
1-butoxybutane                  &-3.88    &2.33      &-1.55     &-0.83      &0.72     \\
trimethoxymethane               &-7.57    &3.51      &-4.06     &-4.42      &-0.36     \\
methoxyethane                   &-4.35    &2.29      &-2.06     &-2.10      &-0.04     \\
1-methoxypropane                &-4.08    &2.24      &-1.84     &-1.66      &0.18     \\
2-methoxypropane                &-4.12    &2.20      &-1.92     &-2.01      &-0.09     \\
1-Ethoxypropane                 &-4.26    &2.32      &-1.94     &-1.81      &0.13     \\
1,3-Dioxolane                   &-6.09    &1.81      &-4.28     &-4.10      &0.18     \\
2,5-dimethyltetrahydrofuran     &-3.86    &1.42      &-2.44     &-2.92      &-0.48     \\
1,1,1-trimethoxyethane          &-7.58    &3.46      &-4.12     &-4.42      &-0.30     \\
2-methoxy-2-methyl-propane      &-3.88    &1.97      &-1.91     &-2.21      &-0.30     \\
1,4-dioxane                     &-7.09    &1.66      &-5.44     &-5.06      &0.38     \\
\hline
RMS                              &    &      &      &      &0.36   \\
\hline
\end{tabular}
\label{ether}
\end{table}

\begin{table}[!ht]
\centering
\caption{The solvation free energy prediction for the alcohol set. All energies are in the unit of kcal/mol.}
\begin{tabular}{rrrrrr}
\cline{1-6}
Name  &$\Delta G^{\rm P}$ &$ G^{\rm NP}$ &${\Delta G}$ &$\Delta G^{\rm Exp}$\cite{Mobley:2014} & Error \\
\hline
ethylene glycol                    &-6.98    &-1.76      &-8.73     &-9.30      &-0.57     \\
butan-1-ol                         &-3.33    &-1.51      &-4.84     &-4.72      &0.12     \\
ethanol                            &-3.49    &-1.47      &-4.96     &-5.00      &-0.04     \\
methanol                           &-3.69    &-1.41      &-5.10     &-5.10      &0.00     \\
propan-1-ol                        &-3.34    &-1.48      &-4.82     &-4.85      &-0.03     \\
propan-2-ol                        &-3.26    &-1.36      &-4.62     &-4.74      & -0.12    \\
pentan-1-ol                        &-3.36    &-1.61      &-4.97     &-4.57      &0.40     \\
2-methylpropan-2-ol                &-3.10    &-1.27      &-4.37     &-4.47      &-0.10     \\
2-methylbutan-2-ol                 &-2.95    &-1.17      &-4.12     &-4.43      &-0.31     \\
2-methylpropan-1-ol                &-3.20    &-1.50      &-4.70     &-4.50      &0.20     \\
butan-2-ol                         &-3.09    &-1.32      &-4.40     &-4.62      &-0.22     \\
cyclopentanol                      &-3.20    &-1.68      &-4.88     &-5.49      &-0.61     \\
4-methylpentan-2-ol                &-2.65    &-1.05      &-3.69     &-3.73      &-0.04     \\
cyclohexanol                       &-3.21    &-1.92      &-5.13     &-5.46      &-0.33     \\
hexan-1-ol                         &-3.43    &-1.53      &-4.96     &-4.40      &0.56     \\
heptan-1-ol                        &-3.48    &-1.62      &-5.09     &-4.21      &0.88     \\
2-methylbutan-1-ol                 &-3.27    &-1.29      &-4.56     &-4.42      &0.14     \\
cycloheptanol                      &-3.07    &-1.89      &-4.96     &-5.48      &-0.52     \\
2-methylpentan-3-ol                &-2.86    &-0.93      &-3.78     &-3.88      &-0.10     \\
pentan-3-ol                        &-3.01    &-1.08      &-4.10     &-4.35      &-0.25     \\
4-Heptanol                         &-2.90    &-1.10      &-3.99     &-4.01      &-0.02     \\
2-methylpentan-2-ol                &-2.93    &-1.08      &-4.00     &-3.92      &0.08     \\
2,3-Dimethyl-2-butanol             &-2.89    &-0.93      &-3.82     &-3.91      &-0.09     \\
hexan-3-ol                         &-3.04    &-1.27      &-4.31     &-4.06      &0.25     \\
pentan-2-ol                        &-3.10    &-1.23      &-4.33     &-4.39      &-0.06     \\
\hline
RMS                              &    &      &      &      &0.33   \\
\hline
\end{tabular}
\label{alcohol}
\end{table}

\begin{table}[!ht]
\centering
\caption{The solvation free energy prediction for the phenol set. All energies are in the unit of kcal/mol.}
\begin{tabular}{rrrrrr}
\cline{1-6}
Name  &$\Delta G^{\rm P}$ &$ G^{\rm NP}$ &${\Delta G}$ &$\Delta G^{\rm Exp}$\cite{Mobley:2014} & Error \\
\hline
3-hydroxybenzaldehyde        &-9.17    &0.39      &-8.78     &-9.52      & -0.74    \\
4-hydroxybenzaldehyde        &-9.60    &0.19      &-9.41     &-8.83      & 0.58    \\
o-cresol                     &-5.32    &-1.04      &-6.36     &-5.90      &0.46     \\
m-cresol                     &-5.71    &-0.86      &-6.57     &-5.49      &1.08     \\
phenol                       &-5.81    &-0.14      &-6.95     &-6.61      &0.34     \\
p-cresol                     &-5.80    &-1.05      &-6.85     &-6.13      &0.72     \\
naphthalen-1-ol              &-5.50    &-0.75      &-6.25     &-7.67      &-1.42     \\
3,4-dimethylphenol           &-5.72    &-0.49      &-6.21     &-6.50      &-0.29     \\
2,5-dimethylphenol           &-5.34    &-0.48      &-5.82     &-5.91      &-0.09     \\
4-tert-butylphenol           &-5.55    &0.86      &-4.69     &-5.91      & -1.22    \\
2,4-dimethylphenol           &-5.55    &-1.03      &-6.58     &-6.01      &0.57     \\
3,5-dimethylphenol           &-5.69    &-0.41      &-6.10     &-6.27      &-0.17     \\
naphthalen-2-ol              &-5.85    &-0.72      &-6.57     &-8.11      &-1.54     \\
2,3-dimethylphenol           &-5.47    &-1.13      &-6.60     &-6.16      &0.44     \\
2,6-dimethylphenol           &-5.07    &-1.07      &-6.14     &-5.26      &0.88     \\
3-ethylphenol                &-5.67    &-0.37      &-6.04     &-6.25      &-0.21     \\
4-propylphenol               &-5.79    &-0.05      &-5.84     &-5.21      &0.63     \\
4-ethylphenol                &-5.76    &-0.48      &-6.24     &-6.13      &0.11     \\
\hline
RMS                              &    &      &      &      &0.76   \\
\hline
\end{tabular}
\label{phenol}
\end{table}
It is interesting to note that for both alkane and alkene molecules, the polar solvation free energy contribution is very small and the nonpolar part dominates the solvation free energy contribution, which explains  why the DG based nonpolar solvation model works extremely well for the solvation free energy prediction of alkane and alkene molecules \cite{ZhanChen:2012}. Further, note that for almost all the alkane molecules, the polar solvation free energies $\Delta G^{\rm P}$ are of magnitude 0.01 kcal/mol, while  alkene molecules have slightly larger magnitude polar free energies, which further verifies that  alkene molecules has a stronger polarity than alkane molecules in general.

Finally, we analyze three classes of polar solute molecules, namely, ether, alcohol, and phenol molecules. Figures  \ref{pred-exp-ether}, \ref{pred-exp-alcohol} and  \ref{pred-exp-phenol}  illustrate  the predicted and experimental solvation free energies for 17 ether, 25 alcohol, and 18 phenol molecules, respectively. Tables  \ref{ether}, \ref{alcohol}  and \ref{phenol}  list the polar, nonpolar, total and experimental solvation free energies for the corresponding families of solute molecules. The RMS errors of these three families are 0.36, 0.33, and 0.76 kcal/mol, respectively.

From the results listed in Tables  \ref{ether}, \ref{alcohol} and \ref{phenol}  we note that for ether molecules, all the nonpolar energies are positive which neutralizes some polar contributions  to the total solvation free energies. For the alcohol molecules, the nonpolar energies are all negative, which enhance the contributions of the polar contributions to the total solvation free energies. Since the surface part is always positive and the volume part is mostly positive, the attractive van der Waals interactions between alcohol molecules and water solvent must be very strong, which explains  that alcohol molecules are easily solvated. As for the phenol molecules, there is a mixed pattern for the nonpolar contributions.

The above study of a large variety of molecules indicates that our DG based solvation model together with the proposed parameter optimization algorithms can provide very accurate predictions of solvation free energies for both polar and nonpolar solute molecules.

%
%
%

\subsection{Five-fold cross validation}

\begin{table}[!ht]
\centering
\caption{The partition of the molecules into sub-groups.}
\begin{tabular}{llllll}
\cline{1-6}
Molecule &Group 1  &Group 2  &Group 3    &Group 4 &Group 5\\
\hline
Alkane   &8    &8       &8 &7   &7 \\
Alkene   &5    &5       &5 &4   &4 \\
Ether    &4    &4       &3 &3   &3 \\
Alcohol  &5    &5       &5 &5   &5 \\
Phenol   &4    &4       &4 &3   &3 \\
\hline
\end{tabular}
\label{mini-group}
\end{table}

Having verified that our DG based solvation model with the optimized parameters   provides very good regression results, we perform a  five-fold cross validation to further illustrate the predictive power of the present method  for independent data sets. Specifically,  the parameters learned from a group of molecules  can be employed for the blind prediction of other molecules.
\begin{figure}[!ht]
\small
\centering
\includegraphics[width=8cm,height=6cm]{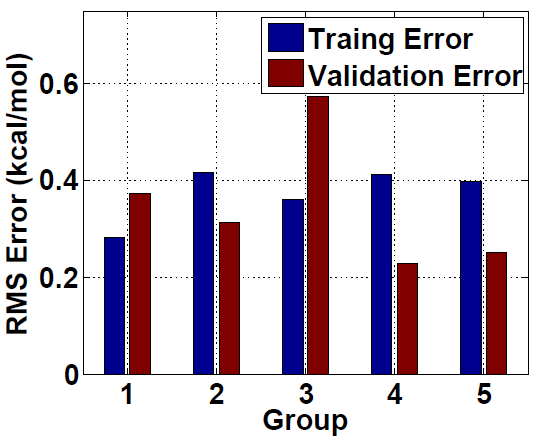}
\caption{The bar plot of the training and validation errors of  alkanes.}
\label{CV_Alkanes}
\end{figure}

\begin{figure}[!ht]
\small
\centering
\includegraphics[width=8cm,height=6cm]{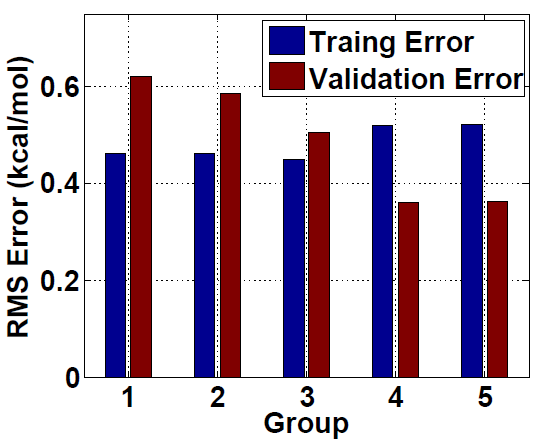}
\caption{The bar plot of the training and validation errors of alkenes.}
\label{CV_Alkenes}
\end{figure}

\begin{figure}[!ht]
\small
\centering
\includegraphics[width=8cm,height=6cm]{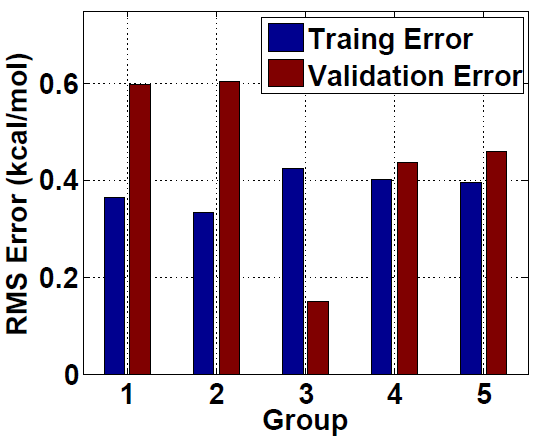}
\caption{The bar plot of the training and validation errors of the ethers.}
\label{CV_ether}
\end{figure}

\begin{figure}[!ht]
\small
\centering
\includegraphics[width=8cm,height=6cm]{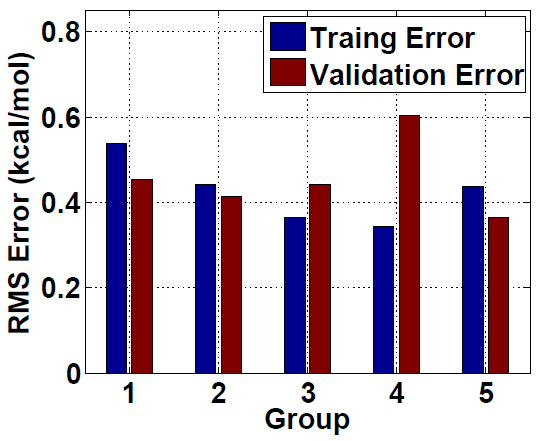}
\caption{The bar plot of the training and validation errors of alcohols.}
\label{CV_Alcohol}
\end{figure}

\begin{figure}[!ht]
\small
\centering
\includegraphics[width=8cm,height=6cm]{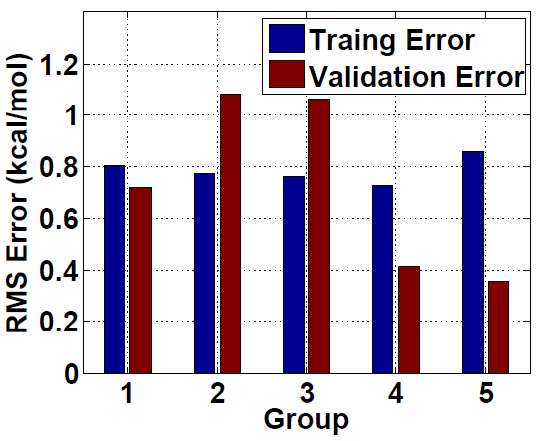}
\caption{The bar plot of the training and validation errors  of  phenols.}
\label{CV_Phenol}
\end{figure}
To perform the five-fold cross validation, each type of molecules is subdivided into five sub-groups as uniformly as possible, Table  \ref{mini-group}   lists the number of molecules in each sub-group for each type of molecules. In our parameters optimization, we leave out one sub-group of molecules and use the rest of molecules to establish our DG based solvation model. The  optimized parameters are then employed for the blind prediction of solvation free energies  of the left out sub-group of molecules.

Figures  \ref{CV_Alkanes}, \ref{CV_Alkenes},\ref{CV_ether}, \ref{CV_Alcohol}, and \ref{CV_Phenol}  demonstrate the cross validation results of the alkane, alkene, ether, alcohol, and phenol molecules, respectively. It is seen that  training and validation errors are similar to each other, which verifies the ability of our model in the blind prediction of solvation free energies.

In the real prediction of the solvation free energy for a given molecule of unknown category, we can first assign it to a given group, and then employ the DG based solvation model with the optimal parameters learned for this specific group for a blind prediction.

\section{Conclusion} \label{Conclusion}

Differential geometry (DG) based solvation models have had a considerable success in solvation analysis \cite{ Wei:1999, ZhanChen:2010a,ZhanChen:2010b, ZhanChen:2011a}.  Particularly, our DG based nonpolar solvation model was shown to offer some of the most accurate  solvation energy predictions of various nonpolar molecules \cite{ZhanChen:2012}. However, our DG based full solvation model  is subject to numerical instability in solving the generalized Laplace-Beltrami (GLB) equation, due to its  coupling with the generalized Poisson Boltzmann (GPB) equation. To stabilize the coupled GLB and GPB equations, a strong constraint on the van der Waals interaction was applied in our earlier work  \cite{ZhanChen:2010a,ZhanChen:2010b, ZhanChen:2011a}, which hinders the parameter optimization of our DG based solvation model. In the present work, we resolve this problem by introducing new parameter optimization algorithms, namely perturbation method and convex optimization,  for the DG based solvation model. New stability conditions are explicitly imposed to the parameter selection,  which guarantees the stability and robustness of solving the GLB equation and leads to constrained optimization of our DG based solvation model. The new optimization algorithms are intensively validated by using  a large number of test molecules, including the SAMPL0 test set \cite{Nicholls:2008},  alkane, alkene, ether, alcohol and phenol types of solutes. Regression results based on our new algorithms are consistent extremely well with experimental data. Additionally, a five-fold cross validation technique is employed to explore  the ability of  our DG based solvation models for the blind prediction of the solvation free energies for a variety of solute molecules.  It is found that the same level of errors is found in the training and validation sets, which confirms our model's predictive power in solvation free energy analysis. The present DG based full solvation model provides a unified framework for analyzing both polar and nonploar molecules.   In our future work, we will develop machine learning approaches for the robust classification of solute molecules of interest into appropriate categories so as to better predict their solvation free energies.

\section*{Acknowledgments}

This work was supported in part by NSF grants IIS-1302285  and DMS-1160352, NIH Grant R01GM-090208, and  MSU Center for Mathematical Molecular Biosciences Initiative. The authors thank Nathan Baker for valuable comments.

\end{document}